\documentclass[reqno]{amsart}
\pdfoutput=1
\usepackage{amssymb}
\usepackage{graphicx}
\usepackage[all]{xy}

\begin{document}
\newcommand*{\threesim}{\mathrel{\vcenter{\offinterlineskip\hbox{$\sim$}\vskip-.35ex\hbox{$\sim$}\vskip-.35ex\hbox{$\sim$}}}}
\newtheorem{prop}{Proposition}[section]
\newtheorem{lemma}[prop]{Lemma}
\newtheorem*{thmA}{Theorem A}
\newtheorem*{thmB}{Theorem B}
\newtheorem{cor}[prop]{Corollary}
\theoremstyle{definition}
\newtheorem{exa}{Example}
\newtheorem{cla}{Claim}
\newtheorem{que}{Question}
\newtheorem*{rmk}{Remark}
\newtheorem*{rmks}{Remarks}
\theoremstyle{remark}
\newtheorem*{ack}{Acknowledgment}
\newtheorem*{pf}{Proof}
\numberwithin{equation}{section}
\title{On Conjugacy Invariants of $D_{\infty}$-Topological Markov Chains}
\author{Sieye Ryu}
\address{Department of Mathematics Ben Gurion University of the Negev,
P.O.B. 653, Be'er Sheva, 8410501, Israel}
\email{Sieye Ryu <sieye@math.bgu.ac.il>}

\begin{abstract}
A $D_{\infty}$-topological Markov chain can be represented by a pair of zero-one square matrices, which is called a flip pair. 
We introduce the concepts of $D_{\infty}$-strong shift equivalence and $D_{\infty}$-shift equivalence, which are equivalence relations between flip pairs.
We investigate the relationships between the existence of a $D_{\infty}$-conjugacy, the existence of a $D_{\infty}$-strong shift equivalence, the existence of a $D_{\infty}$-shift equivalence and the coincidence of the Lind zeta functions. 
\end{abstract}
\maketitle
\section{Introduction}
\label{sec:first}

Time-reversal symmetry arises in many physically-motivated dynamical systems such as classical mechanics, thermodynamics and quantum mechanics. In the early days of dynamical systems, the importance of time-reversal symmetry was recognized in Birkhoff's study of the restricted three body problem \cite{Bi}. After Birkhoff's work, the study of time-reversal symmetry has been extended. For instance, see \cite{D, 1S}. More information on the topic of time-reversal symmetry is provided in the survey \cite{LR}. 

In this paper, we are interested in involutory reversing symmetries defined on  topological Markov chains. 
If an invertible dynamical system $(X, T)$ possesses an involutory reversing symmetries, then it can be regarded as an action of the infinite dihedral group on $X$. 
We study conjugacy invariants of the infinite dihedral group action on topological Markov chains.
For conjugacy invariants of actions of other groups on topological Markov chains (or shifts of finite type), see \cite{BS, BS1}. 

We begin by defining reversals. 
Let $(X, T)$ be an invertible topological dynamical system. 
A homeomorphism $\varphi: X \rightarrow X$ is said to be a \textit{reversal for $(X, T)$} (or \textit{reversing symmetry of $T$}) if $\varphi$ is a conjugacy from $(X, T)$ to $(X, T^{-1})$:
\begin{equation}
\label{eq: 1.1}
\varphi \circ T = T^{-1} \circ \varphi.
\end{equation}
A reversal $\varphi$ for $(X, T)$ is said to be a \textit{flip for $(X, T)$} (or \textit{involutory reversing symmetry of $T$}) if the composition of $\varphi$ with itself is equal to the identity map of $X$:
\begin{equation}
\label{eq: 1.2}
\varphi^2 = \text{id}_X.
\end{equation}
If $\varphi$ is a flip for $(X, T)$, then we call the triple $(X, T, \varphi)$ a \textit{flip system}. 

Suppose that $D_{\infty}$ is the infinite dihedral group generated by $a$ and $b$:
$$D_{\infty}=\langle a, b : ab=ba^{-1}, \, b^2=1 \rangle.$$
To each flip $\varphi$ for $(X, T)$, there corresponds a unique $D_{\infty}$-action $\alpha$ on $X$ such that
$$\alpha(a, x) = T(x) \qquad \text{and} \qquad \alpha(b, x) = \varphi(x).$$
Two flip systems $(X, T, \varphi)$ and $(Y, S, \psi)$ are said to be \textit{$D_{\infty}$-conjugate} if there is a homeomorphism $\theta: X \rightarrow Y$ such that
\begin{equation}
\label{eq: 1.3}
\theta \circ T = S \circ \theta \qquad \text{and} \qquad \theta \circ \varphi = \psi \circ \theta
\end{equation} 
and we write $(X, T, \varphi) \cong (Y, S, \psi)$. 
The homeomorphism $\theta$ is called a \textit{$D_{\infty}$-conjugacy} from $(X, T, \varphi)$ to $(Y, S, \psi)$. 

By (\ref{eq: 1.1}) and (\ref{eq: 1.2}), if $(X, T, \varphi)$ is a flip system, then so are $(X, T, T^n\circ \varphi)$ for all integers $n$. By (\ref{eq: 1.3}), $T^n$ are conjugacies from $(X, T, \varphi)$ to $(X, T, T^{2n}\circ \varphi)$ for all integers $n$.
As a result, if two flip systems $(X, T, \varphi)$ and $(Y, S, S^{2n}\circ\psi)$ are $D_{\infty}$-conjugate for some integer $n$, then $(X, T, \varphi)$ and $(Y, S, \psi)$ are $D_{\infty}$-conjugate. On the other hand, if $(X, T, \varphi)$ and $(Y, S, {S}^{2n+1}\circ\psi)$ are $D_{\infty}$-conjugate for some integer $n$, then $(X, T, \varphi)$ and $(Y, S, S\circ\psi)$ are $D_{\infty}$-conjugate. If two flip systems $(X, T, \varphi)$ and $(Y, S, S \circ \psi)$ are $D_{\infty}$-conjugate, then we say that $(X, T, \varphi)$ and $(Y, S, \psi)$ are \textit{skew $D_{\infty}$-conjugate} and a $D_{\infty}$-conjugacy from $(X, T, \varphi)$ to $(Y, S, S \circ \psi)$ will be called a \textit{skew $D_{\infty}$-conjugacy}. In general, $(Y, S, \psi)$ is not $D_{\infty}$-conjugate to $(Y, S, S\circ \psi)$. (See Proposition \ref{prop: 4.1}.)

Let $\mathcal{A}$ be a finite set. 
If $A$ is a zero-one $\mathcal{A} \times \mathcal{A}$ matrix, then $\textsf{X}_A$ will denote the topological Markov chain (TMC) determined by $A$:
$$\textsf{X}_A = \{ x \in \mathcal{A}^{\mathbb{Z}} : \forall \, i \in \mathbb{Z} \quad A(x_i, x_{i+1}) = 1 \}.$$
In this case, we denote the restriction of the shift map of $\mathcal{A}^{\mathbb{Z}}$ to $\textsf{X}_A$ by $\sigma_A$.
If $\varphi$ is a flip for a TMC $(\textsf{X}_A, \sigma_A)$, then the flip system $(\textsf{X}_A, \sigma_A, \varphi)$ will be called a \textit{$D_{\infty}$-topological Markov chain} or a \textit{$D_{\infty}$-TMC} for short.

In this paper, we will study $D_{\infty}$-strong shift equivalence, $D_{\infty}$-shift equivalence and the Lind zeta functions of $D_{\infty}$-TMCs and investigate the relationships between them. 
We first recall the concepts of strong shift equivalence (SSE) and shift equivalence (SE), which are introduced in \cite{W}. 
SSE and SE are equivalence relations between matrices.
In \cite{W}, Williams proved that two TMCs $(\textsf{X}_A, \sigma_A)$ and $(\textsf{X}_B, \sigma_B)$ are conjugate if and only if there is a SSE between $A$ and $B$. 
In the same paper, he conjectured the existence of a SE guarantees the existence of a SSE.
His conjecture, however, turned out to be false in \cite{KR1, KR2}.

The Artin-Mazur zeta function is also well-known conjugacy invariant of TMCs. 
The formula for the Artin-Mazur zeta function is found in \cite{AM}.
It is well known \cite{LM} that the coincidence of the Artin-Mazur zeta functions of topological Markov chains $(\textsf{X}_A, \sigma_A)$ and $(\textsf{X}_B, \sigma_B)$ does not ensure whether or not $A$ and $B$ are shift equivalent. 
For more details, see Chapter 7 of \cite{LM}.
To illustrate the relationships between these invariants,
we introduce some notation.
Suppose that $A$ and $B$ are zero-one square matrices.
By $A \approx B$, we mean that there is a SSE between $A$ and $B$. 
By $A \sim B$, we mean that there is a SE between them. 
The Artin-Mazur zeta function of $(\textsf{X}_A, \sigma_A)$ is denoted by $\zeta_A(t)$.
The following diagram indicates the relationships between the invariants: 
{\small{$$(\textsf{X}_A, \sigma_A) \cong (\textsf{X}_B, \sigma_B) \quad \Leftrightarrow \quad A \approx B \quad \Rightarrow \quad A \sim B \quad \Rightarrow \quad \zeta_{\sigma_A}(t)=\zeta_{\sigma_B}(t).$$}}

In \cite{KLP}, it is shown that a $D_{\infty}$-TMC is uniquely determined up to $D_{\infty}$-conjugacy by a certain pair $(A, J)$  of zero-one square matrices, which will be called a \textit{flip pair}. (We will give the definition of a flip pair in Section 2.)
We will denote the $D_{\infty}$-TMC determined by a flip pair $(A, J)$ by $(\textsf{X}_A, \sigma_A, \varphi_{A, J})$.
$D_{\infty}$-strong shift equivalence ($D_{\infty}$-SSE) and $D_{\infty}$-shift equivalence ($D_{\infty}$-SE) are equivalence relations between flip pairs, which are analogous to SSE and SE, respectively. 
In Section 2, the notion of half elementary equivalence will be introduced.
If there is a half elementary equivalence between two flip pairs $(A, J)$ and $(B, K)$, then there is a skew $D_{\infty}$-conjugacy (rather than a  $D_{\infty}$-conjugacy) between the corresponding $D_{\infty}$-TMCs $(\textsf{X}_A, \sigma_A, \varphi_{A, J})$ and $(\textsf{X}_B, \sigma_B, \varphi_{B, K})$. 
A $D_{\infty}$-SSE of lag $l$ is a sequence of $l$ half elementary equivalences.
In the same section, we will prove the following theorem, which is analogous to Williams' decomposition theorem:
\begin{thmA}
Suppose that $(A, J)$ and $(B, K)$ are flip pairs.
\newline
\emph{(a)} Two $D_{\infty}$-TMCs $(\textsf{X}_A, \sigma_A, \varphi_{A, J})$ and $(\textsf{X}_B, \sigma_B, \varphi_{B, K})$ are $D_{\infty}$-conjugate if and only if there is a $D_{\infty}$-SSE of lag $2l$ between $(A, J)$ and $(B, K)$ for some positive integer $l$.
\newline
\emph{(b)} Two $D_{\infty}$-TMCs $(\textsf{X}_A, \sigma_A, \varphi_{A, J})$ and $(\textsf{X}_B, \sigma_B, \varphi_{B, K})$ are skew $D_{\infty}$-conjugate if and only if there is a $D_{\infty}$-SSE of lag $2l-1$ between $(A, J)$ and $(B, K)$ for some positive integer $l$.
\end{thmA}

In Section 3, we introduce the notion of $D_{\infty}$-SE.
As in the case of TMCs, if there is a $D_{\infty}$-SSE between $(A, J)$ and $(B, K)$, then there is a $D_{\infty}$-SE between them.

In \cite{L}, Lind introduced a generalization of the Artin-Mazur zeta function, which will be called the Lind zeta function. 
In \cite{KLP}, the Lind zeta functions for $D_{\infty}$-TMCs are expressed in terms of matrices from flip pairs. 
From the formula of it, it will be obvious that the Lind zeta function is a conjugacy invariant of $D_{\infty}$-TMCs.
In Section 4, we will see that not only the existence of a $D_{\infty}$-conjugacy but also the existence of a skew $D_{\infty}$-conjugacy guarantees the coincidence of the Lind zeta functions.

In Section 5, we will show that neither the existence of a $D_{\infty}$-SE nor the coincidence of the Lind zeta functions is a complete conjugacy invariant by proving the following: 
\begin{thmB}
Suppose that $(A, J)$ and $(B, K)$ are flip pairs.
\newline
\emph{(a)} The existence of a $D_{\infty}$-SE between flip pairs $(A, J)$ and $(B, K)$ does not imply that the corresponding $D_{\infty}$-TMCs $(\textsf{X}_A, \sigma_A, \varphi_{A, J})$ and $(\textsf{X}_B, \sigma_B, \varphi_{B, K})$ have the same Lind zeta functions.
\newline
\emph{(b)} The coincidence of the Lind zeta functions of $(\textsf{X}_A, \sigma_A, \varphi_{A, J})$ and $(\textsf{X}_B, \sigma_B, \varphi_{B, K})$ does not guarantee the existence of a $D_{\infty}$-SE from $(A, J)$ to $(B, K)$.
\end{thmB}
We introduce some notation to illustrate the relationships between the existence of a $D_{\infty}$-conjugacy, the existence of a $D_{\infty}$-SSE, the existence of a $D_{\infty}$-SE and the coincidence of the Lind zeta functions.
By $(A, J) \approx (B, K)$, we mean that there is a $D_{\infty}$-SSE between two flip paris $(A, J)$ and $(B, K)$. 
By $(A, J) \sim (B, K)$, we mean that there is a $D_{\infty}$-SE between them.
Throughout the paper, we will assume that $\delta$ is either $0$ or $1$ and $\zeta_{A, J}(t)$ will denote the Lind zeta function for $(\textsf{X}_A, \sigma_A, \varphi_{A, J})$. 

The following diagram indicates the results of this paper:
{\footnotesize{$$\xymatrix
{&{(\textsf{X}_A, \sigma_A, \varphi_{A, J}) \cong (\textsf{X}_B, \sigma_B, \sigma_B^{\delta} \circ \varphi_{B, K})} \ar@2{<->}[d] &\\
&{(A, J) \approx (B, K)} \ar@2{->}[dr]  & \\   
{(A, J) \sim (B, K)}  \ar@2{<-}[ur] & & {\zeta_{A, J}(t) = \zeta_{B, K}(t)}}$$}}
The relationships of the invariants of $D_{\infty}$-TMCs are distinguished from the relationships of those for TMCs.

The outline of this paper is as follows.
In Section 2, we introduce the notions of half elementary equivalence and $D_{\infty}$-SSE. In the same section, we prove Theorem A.
We discuss the notion of $D_{\infty}$-SE in Section 3. 
We devote Section 4 to the Lind zeta functions of $D_{\infty}$-TMCs. 
In Section 5, we prove Theorem B.
In Section 6, we consider a further question.

\begin{ack}
{\rm This research was inspired by \cite{JL} and the author obtained some partial results of this research during the author's doctoral course under the supervision of Young-One Kim. The author deeply thanks Young-One Kim and Jungseob Lee. This research was partially supported by the Israel Science Foundation (grant no. 626/14) and the People Programme (Marie Curie Actions) of the European Union's Seventh Framework Programme (FP7/2007-2013) under
REA grant agreement no. 333598.}
\end{ack}

\section{$D_{\infty}$-Strong Shift Equivalence}
\label{sec:second}

Let $\mathcal{A}$ be a finite set. 
Suppose that $A$ and $J$ are zero-one $\mathcal{A} \times \mathcal{A}$ matrices satisfying
\begin{equation}
\label{eq: 1.4}
JA = A^{\textsf{T}}J \quad \text{and} \quad J^2=I.
\end{equation}
Since $J$ is a zero-one non-singular matrix, it follows that there is a unique map $\tau_J : \mathcal{A} \rightarrow \mathcal{A}$ such that for $a, b \in \mathcal{A}$,
\begin{equation}
\label{eq: 1.5}
J(a, b) = 1 \quad \text{if and only if} \quad \tau_J(a)=b.
\end{equation}
From $J^2=I$, we see that $\tau_J^2 = \text{id}_{\mathcal{A}}$.
Since $AJ=JA^{\textsf{T}}$, we have
\begin{equation}
\label{eq: 1.6}
A(a, b) = A(\tau_J(b), \tau_J(a)) \qquad (a, b \in \mathcal{A}).
\end{equation}
If we define $\varphi_J : \mathcal{A}^{\mathbb{Z}} \rightarrow \mathcal{A}^{\mathbb{Z}}$ by
$$\varphi_J(x)_i = \tau_J(x_{-i}) \qquad (x \in \mathcal{A}^{\mathbb{Z}}; \; i \in \mathbb{Z}),$$
then $\varphi_J$ is a flip for $(\mathcal{A}^{\mathbb{Z}}, \sigma)$. 
By (\ref{eq: 1.6}), we have $\varphi_J(\textsf{X}_A) = \textsf{X}_A$.
Thus, the restriction of $\varphi_J$ to $\textsf{X}_A$ also becomes a flip for $(\textsf{X}_A, \sigma_A)$.
We denote the restriction by $\varphi_{A, J}$. 
A pair $(A, J)$ of zero-one square matrices satisfying (\ref{eq: 1.4}) will be called a \textit{flip pair}. 
 
A flip $\varphi$ for a shift space $(X, \sigma_X)$ is said to be a \textit{one-block flip} if
$$x, x' \in X \;\; \text{and} \;\; x_0=x'_0 \quad \Rightarrow \quad \varphi(x)_0=\varphi(x')_0.$$
In this case, there is a unique map $\tau:\mathcal{A} \rightarrow \mathcal{A}$ such that 
$\tau ^2 = \text{id}_{\mathcal{A}}$ and that
$$\varphi(x)_i=\tau(x_{-i}) \qquad (x \in X ; \; i \in \mathbb{Z}).$$
We call $\tau$ the \textit{symbol map} of $\varphi$.
If $(A, J)$ is a flip pair, then it is obvious that $\varphi_{A, J}$ is a one-block flip for $(\textsf{X}_A, \sigma_A)$ whose symbol map is $\tau_J$ satisfying (\ref{eq: 1.5}).

The following lemma is proved in \cite{KLP, R}.

\begin{lemma}
\label{lemma: 2.1}
If $(X, \sigma_X, \varphi)$ is a $D_{\infty}$-TMC, then there is a flip pair $(A, J)$ such that $(X, \sigma_X, \varphi) \cong (\textsf{X}_A, \sigma_A, \varphi_{A, J})$.
\end{lemma}

Let $(A, J)$ and $(B, K)$ be flip pairs. A pair $(R, S)$ of zero-one matrices satisfying
$$A=RS, \quad B=SR, \quad \text{and} \quad S=KR^{\textsf{T}}J$$
is said to be a \textit{half elementary equivalence from $(A, J)$ to $(B, K)$}. 
We note that $S = KR^{\textsf{T}}J$ is equivalent to $R=JS^{\textsf{T}}K$.
If there is a half elementary equivalence from $(A, J)$ to $(B, K)$, then we write $(R, S) : (A, J) \threesim (B, K)$.

\begin{prop}
\label{prop: 2.2}
If $(R, S) : (A, J) \threesim (B, K)$, then $(\textsf{X}_A, \sigma_A, \varphi_{J, A})$ is skew-$D_{\infty}$ conjugate to $(\textsf{X}_B, \sigma_B, \varphi_{K, B})$, that is, $(\textsf{X}_A, \sigma_A, \varphi_{J, A}) \cong (\textsf{X}_B, \sigma_B, \sigma_B\circ \varphi_{K, B})$. 
\end{prop}

Before proving Proposition \ref{prop: 2.2}, we introduce some notation.
Suppose that $X$ is a shift space over $\mathcal{A}$. 
If $n$ is a positive integer, we denote the set of all admissible $n$-blocks of $X$ by $\mathcal{B}_n(X)$: 
$$\mathcal{B}_n(X) = \{w \in \mathcal{A}^n : \; w \text{ occurs in } x \text{ for some } x \in X\}.$$

\begin{pf}
Since $R$ and $S$ are zero-one and $A=RS$, it follows that for all $a_1a_2 \in \mathcal{B}_2(\textsf{X}_A)$, there is a unique $b\in \mathcal{B}_1(\textsf{X}_B)$ such that 
$$R(a_1, b) = S(b, a_2)=1.$$
We define the block map $\Gamma_{R, S} : \mathcal{B}_2(\textsf{X}_A) \rightarrow \mathcal{B}_1(\textsf{X}_B)$ to be
$$\Gamma_{R, S}(a_1a_2) = b \quad \Leftrightarrow \quad R(a_1, b) = S(b, a_2)=1$$
for $a_1 a_2 \in \mathcal{B}_2(\textsf{X}_A)$ and $b \in \mathcal{B}_1(\textsf{X}_B)$.
If we define the map $\gamma_{R, S} : (\textsf{X}_A, \sigma_A) \rightarrow (\textsf{X}_B, \sigma_B)$ by
$$\gamma_{R, S}(x)_i = \Gamma_{R, S} \left(x_i x_{i+1} \right) \qquad (x \in \textsf{X}_A; \; i \in \mathbb{Z}),$$
then we have $\gamma_{R, S} \circ \sigma_A = \sigma_B \circ \gamma_{R, S}$.

Since $(S, R): (B, K) \threesim (A, J)$, the block map $\Gamma_{S, R} : \mathcal{B}_2(\textsf{X}_B) \rightarrow \mathcal{B}_1(\textsf{X}_A)$ and the map $\gamma_{S, R}: (\textsf{X}_B, \sigma_B) \rightarrow (\textsf{X}_A, \sigma_A)$ are also well-defined:
$$\Gamma_{S, R}(b_1b_2) = a \quad \Leftrightarrow \quad S(b_1, a) = R(a, b_2)=1 \quad (b_1b_2 \in \mathcal{B}_2(\textsf{X}_B), a \in \mathcal{B}_1(\textsf{X}_A))$$
and
$$\gamma_{S, R}(y)_i = \Gamma_{S, R} \left(y_i y_{i+1} \right) \qquad (y \in \textsf{X}_B; \; i \in \mathbb{Z}).$$
Since $\gamma_{S, R} \circ \gamma_{R, S} = \sigma_A$ and $\gamma_{R, S} \circ \gamma_{S, R} = \sigma_B$, it follows that $\gamma_{R, S}$ is one-to-one and onto.

It remains to show that 
\begin{equation}
\label{eq: 2.5}
\gamma_{R, S} \circ \varphi_{J, A} = \left( \sigma_B \circ \varphi_{K, B} \right) \circ \gamma_{R, S}.
\end{equation}
Since $S = KR^{\textsf{T}}J$, it follows that 
$$S(b, a) = 1 \quad \Leftrightarrow \quad R(\tau_J(a), \tau_K(b))=1 \qquad (a \in \mathcal{B}_1(\textsf{X}_A),  b\in \mathcal{B}_1(\textsf{X}_B))$$
or equivalently,
$$R(a, b) = 1 \quad \Leftrightarrow \quad S(\tau_K(b), \tau_J(a))=1 \qquad (a \in \mathcal{B}_1(\textsf{X}_A),  b\in \mathcal{B}_1(\textsf{X}_B)).$$
Thus, we obtain
\begin{equation}
\label{eq: 2.6}
\Gamma_{R, S}(a_1 a_2) =b \quad \Leftrightarrow \quad \Gamma_{R, S} \left( \tau_J(a_2) \tau_J(a_1) \right) = \tau_K(b) \qquad \left( a_1 a_2 \in \mathcal{B}_2 (\textsf{X}_A) \right)
\end{equation}
By (\ref{eq: 2.6}), we have
\begin{eqnarray*}
\gamma_{R, S} \circ \varphi_{J, A} (x)_i
&=& \Gamma_{R, S}(\tau_J(x_{-i}) \tau_J(x_{-i-1}))
= \tau_K \Gamma_{R, S}(x_{-i-1}x_{-i})\\
&=& \varphi_{K, B} \circ \gamma_{R, S} (x)_{i+1}
= \left( \sigma_B \circ \varphi_{K, B} \right) \circ \gamma_{R, S} (x)_i 
\end{eqnarray*}
and (\ref{eq: 2.5}) is proved. (See also the diagrams below.)
\hfill $\Box$
\end{pf}

The following diagram indicates (\ref{eq: 2.6}):
\begin{displaymath}
\xymatrix{a_1 \ar[rr]^A \ar[dr]_R &&a_2 &\Leftrightarrow& \tau_J(a_2) \ar[rr]^A \ar[dr]_R &&\tau_J(a_1) \\&b \ar[ur]^S& & & & \tau_K(b) \ar[ur]^S&}                   
\end{displaymath}
If we underline the zero-th coordinate of a bi-infinite sequence $x$, that is,
$$x=\cdots x_{-2} x_{-1} \underline{x_0} x_1 x_2 \cdots,$$
then the following diagram indicates (\ref{eq: 2.5}).
\begin{displaymath}
\xymatrix{\cdots \qquad x_{-1} \ar[rr]^A \ar[dr]_R && \underline{x_0} \ar[rr]^A \ar[dr]_R &&  {x_1} \qquad \cdots \\        \cdots   &    y_{-1} \ar[rr]_B \ar[ur]^S &  &\underline{y_0} \ar[ur]^S &  \cdots}                   
\end{displaymath}  
$$\Updownarrow$$
\begin{displaymath}
\xymatrix{\cdots \qquad {\tau_J(x_{1})} \ar[rr]^A \ar[dr]_R && \underline{\tau_J(x_0)} \ar[rr]^A \ar[dr]_R &&  \tau_J(x_{-1}) \qquad \cdots && \\
        \cdots &    {\tau_K(y_{0})} \ar[rr]_B \ar[ur]^S &  &\underline{\tau_K(y_{-1})}  \ar[ur]^S   &\cdots}                         
\end{displaymath}

Let $(A, J)$ and $(B, K)$ be flip pairs. 
A sequence of $l$ half elementary equivalences 
$$(R_1, S_1) : (A, J) \threesim (A_2, J_2),$$
$$(R_2, S_2) : (A_2, J_2) \threesim (A_3, J_3),$$
$$\vdots$$
$$(R_l, S_l): (A_l, J_l) \threesim (B, K)$$ 
is said to be a  \textit{$D_{\infty}$-SSE of lag $l$ from $(A, J)$ to $(B, K)$}.
If there is a $D_{\infty}$-SSE of lag $l$ from $(A, J)$ to $(B, K)$, then we say that $(A, J)$ is $D_{\infty}$-strong shift equivalent to $(B, K)$ and write $(A, J) \approx (B, K)$ (lag $l$), or simply $(A, J) \approx (B, K)$. 

By Proposition \ref{prop: 2.2}, it is clear that
\begin{equation}
\label{eq: 2.7}
(A, J) \approx (B, K) \; (\text{lag} \; l) \quad \Rightarrow \quad (\textsf{X}_A, \sigma_A, \varphi_{J, A}) \cong (\textsf{X}_B, \sigma_B, {\sigma_B}^l \circ \varphi_{K, B}).
\end{equation}
Because ${\sigma_B}^l$ is a conjugacy from $(\textsf{X}_B, \sigma_B, \varphi_{K, B})$ to $(\textsf{X}_B, \sigma_B, {\sigma_B}^{2l} \circ \varphi_{K, B})$, the implication in (\ref{eq: 2.7}) can be rewritten as follows:
\begin{equation}
\label{eq: 2.8}
(A, J) \approx (B, K) \; (\text{lag} \; 2l) \quad \Rightarrow \quad (\textsf{X}_A, \sigma_A, \varphi_{J, A}) \cong (\textsf{X}_B, \sigma_B, \varphi_{K, B})
\end{equation}
and 
$$(A, J) \approx (B, K) \; (\text{lag} \; 2l-1) \quad \Rightarrow \quad (\textsf{X}_A, \sigma_A, \varphi_{J, A}) \cong (\textsf{X}_B, \sigma_B, \sigma_B \circ \varphi_{K, B}).$$  
We will see that $(\textsf{X}_B, \sigma_B, \varphi_{K, B})$ is not $D_{\infty}$-conjugate to $(\textsf{X}_B, \sigma_B, \sigma_B \circ \varphi_{K, B})$ but they share the same Lind zeta functions in Proposition \ref{prop: 4.1}.

In order to prove Theorem A, we introduce the notion of a higher block shift-flip system.
We first indicate some notation.
Suppose that $\mathcal{A}$ is a finite set.
When $n$ is a positive integer, we define the $n$-th mirror map $\rho_n : \mathcal{A}^n \rightarrow \mathcal{A}^n$, the $n$-th initial map $i_n : \bigcup_{k=n}^{\infty} \mathcal{A}^k \rightarrow \mathcal{A}^n$ and the $n$-th terminal map $t_n : \bigcup_{k=n}^{\infty} \mathcal{A}^k \rightarrow \mathcal{A}^n$ by
$$\rho_n(a_1 a_2 \cdots a_n) = a_n \cdots a_2 a_1 \qquad (a_1 a_2 \cdots a_n \in \mathcal{A}^n),$$
$$i_n(a_1 a_2 \cdots a_m) = a_1 a_2 \cdots a_n \qquad (a_1 a_2 \cdots a_m \in \mathcal{A}^m; \; m \geq n)$$
and
$$t_n(a_1 a_2 \cdots a_m) = a_{m-n+1} a_{m-n+2} \cdots a_m \qquad (a_1 a_2 \cdots a_m \in \mathcal{A}^m; \; m \geq n).$$
For notational convenience, we drop $n$ and denote the $n$-th mirror map by $\rho$. It would be clear what $\rho$ indicates in the context. 
We also denote the restrictions of $\rho, i_n, t_n$ to $\mathcal{B}_n(X)$  by $\rho, i_n, t_n$ when $X$ is a shift space over $\mathcal{A}$.

Suppose that $(X, \sigma_X)$ is a shift space and that $\varphi$ is a one-block flip for $(X, \sigma_X)$ with a symbol map $\tau$.
The $n$-th higher block system $(X_n, \sigma_n)$ has a natural one-block flip $\varphi_n$ whose symbol map is $\rho \circ \tau_n$:
$$\varphi_n(x)_i = (\rho \circ \tau_n)(x_{-i}) \qquad (x \in X_n; \;i \in \mathbb{Z}).$$
For instance, if $x \in X$, then we have
$$\cdots \left[ \begin{array}{c} x_{0} \\ x_{-1} \end{array}\right] \underline{\left[ \begin{array}{c} x_{1} \\ x_{0} \end{array}\right]} \left[ \begin{array}{c} x_{2} \\ x_{1} \end{array}\right] \cdots \in X_2$$
and
$$\varphi_2\left(\cdots \left[ \begin{array}{c} x_{0} \\ x_{-1} \end{array}\right] \underline{\left[ \begin{array}{c} x_{1} \\ x_{0} \end{array}\right]} \left[ \begin{array}{c} x_{2} \\ x_{1} \end{array}\right] \cdots \right) = \cdots \left[ \begin{array}{c} 
\tau(x_{1}) \\ \tau(x_{2}) \end{array}\right] \underline{\left[ \begin{array}{c} \tau(x_{0}) \\ \tau(x_{1}) \end{array}\right]} \left[ \begin{array}{c} \tau(x_{-1}) \\ \tau(x_{0}) \end{array}\right] \cdots.$$
For notational simplicity, we drop the subscript $n$ and write $\tau=\tau_n$.
$$\tau(a_1 a_2 \cdots a_n) = \tau(a_1) \tau(a_2) \cdots \tau(a_n).$$
The shift-flip system $(X_n, \sigma_n, \varphi_n)$ will be called the \textit{$n$-th higher block shift-flip system of $(X, \sigma_X, \varphi)$}.

For each positive integer $n$, we define zero-one  $\mathcal{B}_n(\textsf{X}_A) \times \mathcal{B}_n(\textsf{X}_A)$ matrices $A_n$ and $J_n$ by
$$A_n(u, v) = \begin{cases} 1 \qquad \text{if } t_{n-1}(u)=i_{n-1}(v),\\0 \qquad \text{otherwise}\end{cases} \qquad \big(u, v \in \mathcal{B}_n(\textsf{X}_A)\big)$$
and
$$J_n(u, v) = \begin{cases} 1 \qquad \text{if } v = (\rho \circ 
\tau)(u),\\0 \qquad \text{otherwise}\end{cases} \qquad \big(u, v \in \mathcal{B}_n(\textsf{X}_A)\big)$$
so that
$(A_n, J_n)$ is a flip pair for $(X_n, \sigma_n, \varphi_n)$.

\begin{lemma}
\label{lemma: 2.3}
For each positive integer $n$, we have
$$(A_1, J_1) \approx (A_{n+1}, J_{n+1})\; (\emph{lag} \, n).$$
\end{lemma}

\begin{pf}
For each $k=1, 2, \cdots, n$,
we define a zero-one $\mathcal{B}_k(\textsf{X}_A) \times \mathcal{B}_{k+1}(\textsf{X}_A)$ matrix $R_k$ and a zero-one $ \mathcal{B}_{k+1}(\textsf{X}_A) \times \mathcal{B}_k (\textsf{X}_A)$ matrix $S_k$ by 
$$R_k(u, v)=\begin{cases} 1 \qquad \text{if} \; u=i_{k}(v), \\ 0 \qquad \text{otherwise}, \end{cases} \qquad \big(u \in \mathcal{B}_k(\textsf{X}_A), \, v \in \mathcal{B}_{k+1}(\textsf{X}_A)\big)$$
and
$$S_k(v, u)=\begin{cases} 1 \qquad \text{if} \; u=t_k(v), \\ 0 \qquad \text{otherwise} \end{cases} \qquad \big(u \in \mathcal{B}_k(\textsf{X}_A), \, v \in \mathcal{B}_{k+1}(\textsf{X}_A)\big).$$
It is straightforward to see that $(R_k, S_k): (A_k, J_k) \threesim (A_{k+1}, J_{k+1})$ for each $k$.
\hfill $\Box$
\end{pf}

In the rest of the section, we prove (a) in Theorem A. 
(b) is an immediate consequence of (a) and Lemma \ref{lemma: 2.3}.
We denote the flip pairs for the $n$-th higher block shift-flip systems of $(\textsf{X}_{A}, \sigma_A, \varphi_{A, J})$ and $(\textsf{X}_{B}, \sigma_B, \varphi_{B, K})$ by $(A_n, J_n)$ and $(B_n, K_n)$, respectively.

Suppose that $\psi :(\textsf{X}_{A}, \sigma_A, \varphi_{A, J}) \rightarrow (\textsf{X}_{B}, \sigma_B, \varphi_{B, K})$ is a conjugacy. Then there are nonnegative integers $s$ and $t$ and a block map $\Psi : \mathcal{B}_{s+t+1}(\textsf{X}_A) \rightarrow \mathcal{B}_1(\textsf{X}_B)$ such that
$$\psi(x)_i = \Psi(x_{[i-s, i+t]}) \qquad (x \in \textsf{X}_A; \; i \in \mathbb{Z}).$$
We may assume that $s+t$ is even by extending window size if necessary. 
By Lemma \ref{lemma: 2.3}, there is a $D_{\infty}$-SSE of lag $(s+t)$ from $(A, J)$ to $(A_{s+t+1}, J_{s+t+1})$.
From (\ref{eq: 2.8}), the ($s+t+1$)-th higher block code $h$
is a $D_{\infty}$-conjugacy. 
It is clear that there is a one-block conjugacy $\psi'$ induced by $\psi$ from the ($s+t+1$)-th higher block shift-flip system to $(\textsf{X}_{B}, \sigma_B, \varphi_{B, K})$ satisfying $\psi = \psi' \circ h$. We need to prove that there is a $D_{\infty}$-SSE of lag $2l$ from $(A_{s+t+1}, J_{s+t+1})$ to $(B, K)$ for some positive integer $l$. We may assume that $s=t=0$ and prove that there is a $D_{\infty}$-SSE of lag $2l$ from $(A, J)$ to $(B, K)$ for some positive integer $l$.

If $\psi^{-1}$ is the inverse of $\psi$, there is a nonnegative integer $m$ such that
\begin{equation}
\label{eq: 2.10}
y, y' \in \textsf{X}_B \;\; \text{and} \;\; y_{[-m, m]}=y'_{[-m, m]} \quad \Rightarrow \quad \psi^{-1}(y)_0=\psi^{-1}(y')_0.
\end{equation}
For each $k=1, 2, \cdots, 2m+1$, we define a set $\mathcal{A}_k$ by
$$\mathcal{A}_k = \left\{ \left[\begin{array}{c} v \\ w \\ u \end{array} \right]: u, v \in \mathcal{B}_{i}(\textsf{X}_B), w \in \mathcal{B}_{j}(\textsf{X}_A) \text{ and } u \Psi(w) v \in \mathcal{B}_k(\textsf{X}_{B})\right\},$$
where $i=\lfloor \frac{k-1}{2} \rfloor$ and $j=k-2\lfloor\frac{k-1}{2} \rfloor$.
Next, we define $\mathcal{A}_k \times \mathcal{A}_k$ matrices $M_k$ and $F_k$ to be
\begin{eqnarray*}
 M_k \Bigg( \left[\begin{array}{c} v \\ w \\ u \end{array} \right], \left[\begin{array}{c} v' \\ w' \\ u' \end{array} \right] \Bigg) = 1 \qquad
&\Leftrightarrow& \qquad \left[\begin{array}{c} v \\ \Psi(w) \\ u \end{array} \right] \left[\begin{array}{c} v' \\ \Psi(w') \\ u' \end{array} \right] \in \mathcal{B}_2(\textsf{X}_{B_k}) \\
&&\\
&& \qquad \quad \text{and} \quad ww' \in \mathcal{B}_2(\textsf{X}_{A_j})
\end{eqnarray*}
and
$$
F_k \Bigg( \left[\begin{array}{c} v \\ w \\ u \end{array} \right], \left[\begin{array}{c} v' \\ w' \\ u' \end{array} \right] \Bigg)=1 
\qquad \Leftrightarrow \qquad \begin{array}{c} \\ u' = (\rho \circ \tau_{K}) (v), \; w'= (\rho \circ \tau_{J})(w) \\ \\
\qquad \quad \text{and}\;\; v'= (\rho \circ \tau_{K})(u) \end{array}
$$
for all $$\left[\begin{array}{c} v \\ w \\ u \end{array} \right] \;, \; \left[\begin{array}{c} v' \\ w' \\ u' \end{array} \right] \in \mathcal{A}_k.$$
A direct calculation shows that $(M_k, F_k)$ is a flip pair for each $k$.

Now, we construct a $D_{\infty}$-SSE of lag $2m$ from $(A, J)$ to $(M_{2m+1}, F_{2m+1})$.
Define a zero-one $\mathcal{A}_k \times \mathcal{A}_{k+1}$ matrix $R_k$ and a zero-one $\mathcal{A}_{k+1} \times \mathcal{A}_k$ matrix $S_k$ to be
$$ R_k \Bigg( \left[\begin{array}{c} v \\ w \\ u \end{array} \right], \left[\begin{array}{c} v' \\ w' \\ u' \end{array} \right] \Bigg)  = 1 
\qquad \Leftrightarrow \qquad  \begin{array}{c} \\ u \Psi(w) v = i_k \left( u' \Psi(w') v' \right) \\ \\
\quad \quad \text{and} \;\; t_1(w)=i_1(w') \end{array}$$
and
$$S_k\Bigg( \left[\begin{array}{c} v' \\ w' \\ u' \end{array} \right], \left[\begin{array}{c} v \\ w \\ u \end{array} \right] \Bigg)  =  1 \qquad
\Leftrightarrow \qquad \begin{array}{c} \\ t_k \left( u'\Psi(w')v' \right)=u\Psi(w)v \\ \\
\quad \quad \text{and} \;\; t_1(w')=i_1(w), \end{array}$$
for all 
$$\left[\begin{array}{c} v \\ w \\ u \end{array} \right] \in \mathcal{A}_{k} \qquad \text{and} \qquad \left[\begin{array}{c} v' \\ w' \\ u' \end{array} \right] \in \mathcal{A}_{k+1}.$$
A direct calculation shows that
$$(R_k, S_k):(M_k, F_k)\threesim (M_{k+1}, F_{k+1}).$$
Because $M_1=A$ and $F_1=J$, we obtain
\begin{equation}
\label{eq: 2.11}
(A, J) \approx (M_{2m+1}, F_{2m+1})\; (\text{lag} \;2m).
\end{equation}

Finally, (\ref{eq: 2.10}) implies that the $D_{\infty}$-TMC determined by the flip pair $(M_{2m+1}, F_{2m+1})$ is equal to the $(2m+1)$-th higher block shift-flip system of $(\textsf{X}_B, \sigma_B, \varphi_{K, B})$ by recoding of symbols.
From Lemma \ref{lemma: 2.3}, we have
\begin{equation}
\label{eq: 2.12}
(B, K) \approx (M_{2m+1}, F_{2m+1})\; (\text{lag} \;2m).
\end{equation}
From (\ref{eq: 2.11}) and (\ref{eq: 2.12}), it follows that 
$$(A, J) \approx (B, K)\; (\text{lag} \;4m)$$
\hfill $\Box$

\section{$D_{\infty}$-Shift Equivalence}
\label{sec:third}

Let $(A, J)$ and $(B, K)$ be flip pairs and let $l$ be a positive integer.
A \textit{$D_{\infty}$-SE of lag $l$ from $(A, J)$ to $(B, K)$} is a pair $(R, S)$ of nonnegative integral matrices satisfying
$$A^{l}=RS, \quad B^{l}=SR, \quad AR=RB, \quad \text{and} \quad S=KR^{\textsf{T}}J.$$
We observe that $AR=RA$, $S=KR^{\textsf{T}}J$ and the fact that $(A,J)$ and $(B,K)$ are flip pairs imply $SA=BS$.
If there is a $D_{\infty}$-SE of lag $l$ from $(A, J)$ to $(B, K)$, then we say that $(A, J)$ is $D_{\infty}$-shift equivalent to $(B, K)$ and write $(R, S): (A, J) \sim (B, K)$ (lag $l$), or simply $(A, J) \sim (B, K)$.
Suppose that  
$$(R_1, S_1), (R_2, S_2), \cdots, (R_l, S_l)$$ 
is a $D_{\infty}$-SSE of lag $l$ from $(A, J)$ to $(B, K)$. 
If we set
$$R=R_1 R_2 \cdots R_l \qquad \text{and} \qquad S=S_l \cdots S_2 S_1,$$
then $(R, S)$ is a $D_{\infty}$-SE of lag $l$ from $(A, J)$ to $(B, K)$.
Hence,
$$(A, J) \approx (B, K) \;\; (\text{lag } l) \qquad \Rightarrow \qquad (A, J) \sim (B, K) \; (\text{lag } l).$$

\section{The Lind Zeta Functions}
\label{sec:fourth}

Suppose that $G$ is a group and that $\alpha$ is a $G$-action on $X$. Let $\mathcal{F}$ denote the set of finite index subgroups. For each $H \in \mathcal{F}$, we set
$$p_H(\alpha) = |\{ x \in X : \forall \, h \in H \; \alpha(h, x) = x \}|.$$ 
The Lind zeta function $\zeta_{\alpha}$ of the action $\alpha$ is defined by 
$$\zeta_{\alpha}(t) = \exp \left( \sum_{H \in \mathcal{F}} \frac{p_H(\alpha)}{|G/H|}\, t^{|G/H|}\right).$$
It is clear that if $\alpha : \mathbb{Z} \times X \rightarrow X$ is given by $\alpha(n, x) = T^n(x)$, then the Lind zeta function $\zeta_{\alpha}$ becomes the Artin-Mazur zeta function $\zeta_T$ \cite{AM}. Lind defined this function in \cite{L} for the case $G = \mathbb{Z}^d$.

We briefly discuss an explicit formula for the Lind zeta function of a flip system $(X, T, \varphi)$. For more details, see \cite{KLP}. Every finite index subgroup of $D_{\infty} =\langle a, b : ab=ba^{-1} \;\, \text{and} \;\, b^2=1\rangle$ can be written in one and only one of the following forms: 
$$\langle a^m \rangle \qquad \text{or} \qquad \langle a^m, a^k b \rangle \qquad (m=1, 2, \cdots\,; k=1, 2, \cdots, m-1)$$
and has index
$$|G_2/ \langle a^m \rangle| = 2m \qquad \text{or} \qquad |G_2/ \langle a^m, a^k b \rangle| = m.$$
If $m$ is a positive integer, then the number of periodic points in $X$ of period $m$ will be denoted by $p_m(T)$:
$$p_m(T) = |\{ x \in X : T^m(x) = x \}|.$$
If $m$ is a positive integer and $n$ is an integer, then $p_{m,n}(T, \varphi)$ will denote the number of points in $X$ fixed by $T^m$ and $T^n \circ \varphi$:
$$p_{m, n}(T, \varphi) = |\{ x \in X : T^m(x) = T^n \circ \varphi (x) = x \}|.$$
Now, we obtain
$$\zeta_{T, \varphi}(t) = \exp \Big(\sum_{m=1}^{\infty} \, \frac{p_m(T)}{2m}t^{2m} +\sum_{m=1}^{\infty} \, \sum_{k=0}^{m-1}\, \frac{p_{m,k}(T, \varphi)}{m}t^{m} \Big). $$
The definition of a flip gives that
$$p_{m,n}(T, \varphi) = p_{m,n+m}(T, \varphi) = p_{m,n+2}(T, \varphi)$$
and this implies 
\begin{eqnarray}
&& p_{m,n}(T, \varphi) = p_{m,0}(T, \varphi) \qquad \text{ if } m \text{ is odd}, \label{eq: 4.1}\\ && p_{m,n}(T, \varphi)= p_{m,0}(T, \varphi) \qquad \text{ if } m \text{ and } n \text{ are even}, \nonumber \\ && p_{m,n}(T, \varphi)= p_{m,1}(T, \varphi) \qquad \text{ if } m \text{ is even and } n \text{ is odd}. \nonumber 
\end{eqnarray}
Hence, we obtain
$$\sum_{k=0}^{m-1} \, \frac{p_{m,n}(T, \varphi)}{m} = \begin{cases} \, p_{m,0}{T, F} \qquad \qquad \qquad \qquad \text{if } m \text{ is odd},\\ \\ \, \displaystyle{\frac{p_{m,0}(T, \varphi)+p_{m,1}(T, \varphi)}{2} \qquad \text{if } m \text{ is even}.}\end{cases}$$
The Lind zeta function $\zeta_{T, \varphi}$ of a flip system $(X, T, \varphi)$ is given by
$$\zeta_{\alpha}(t) = {\zeta_T(t^2)}^{1/2} \exp \left( G_{T, \varphi}(t) \right),$$
where
$\zeta_T$ is the Artin-Mazur zeta function of $(X, T)$ and
$$G_{T, \varphi}(t) = \sum_{m=1}^{\infty} \, \left( p_{2m-1, 0}(T, \varphi) \, t^{2m-1} + \frac{p_{2m, 0}(T, \varphi)+p_{2m, 1}(T, \varphi)}{2} \, t^{2m}\right).$$
It is evident if the flip systems $(X, T, \varphi)$ and $(X', T', \varphi ')$ are $D_{\infty}$-conjugate, then 
$$p_m(T) = p_m(T')$$
and
$$p_{m, n}(T, \varphi) = p_{m,n}(T', \varphi ')$$
for all positive integers $m$ and integers $n$. 
As a consequence, the Lind zeta function is a conjugacy invariant.

\begin{prop}
\label{prop: 4.1}
If $(X, T, \varphi)$ is a flip system, then
$$p_{2m-1, 0}(T, \varphi) = p_{2m-1, 0}(T, T \circ \varphi),$$
$$p_{2m, 0}(T, \varphi) = p_{2m, 1}(T, T \circ \varphi), \quad \text{and}$$
$$p_{2m, 1}(T, \varphi) = p_{2m, 0}(T, T \circ \varphi) \quad (m=1, 2, \cdots),$$
As a consequence, the Lind zeta functions for $(X, T, \varphi)$ and $(X, T, T \circ \varphi)$ are the same.
\end{prop}

\begin{pf}
The last equality is trivially true. To prove the first two equalities, we observe that 
$x$ is fixed by $T^m$ and $\varphi$ if and only if
$Tx$ is fixed by $T^m$ and $T \circ (T\circ \varphi)$ for all positive integers $m$:
$$T^{m}(x) = \varphi(x)=x \quad \Leftrightarrow \quad T^{m}(Tx)=T \circ (T \circ \varphi)(Tx) = Tx.$$ 
Thus, we have
\begin{equation}
\label{eq: 4.2}
p_{m, 0}(T, \varphi)=p_{m, 1}(T, T \circ \varphi) \qquad (m=1, 2, \cdots).
\end{equation}
Replacing $m$ with $2m$ yields the second equality.
From (\ref{eq: 4.1}) and (\ref{eq: 4.2}), the first one follows.
\hfill $\Box$
\end{pf}

\begin{cor}
\label{cor: 4.2}
Let $(A, J)$ and $(B, K)$ be flip pairs. If there is a $D_{\infty}$-SSE of lag $l$ from $(A, J)$ to $(B, K)$ for some positive integer $l$, then the Lind zeta functions for  
$(\textsf{X}_A, \sigma_A, \varphi_{J, A})$ and $(\textsf{X}_B, \sigma_B, \varphi_{K, B})$ are the same.
\end{cor}

Let $(A, J)$ be a flip pair again. 
It is well-known \cite{LM} that 
$$p_{m}(\sigma_A) = \text{tr}(A^m) \qquad (m = 1, 2, \cdots).$$
In \cite{KLP}, a similar formula for $p_{m, \delta}(\sigma_X, \varphi)$ is established when $m$ is a positive integer.
In order to present it, we indicate notation. 
If $M$ is a square matrix, then $\Delta_M$ will denote the column vector whose $i$-th coordinates are identical with $i$-th diagonal entries of $M$, that is
$$\Delta_M(i) = M_{ii}.$$
For instance, if $I$ is the $2 \times 2$ identity matrix, then 
$$\Delta_I = \left[\begin{array}{c} 1 \\ 1 \end{array}\right].$$
With this notation, we obtain a formula for $p_{m,n}(\sigma_X, \varphi)$ as follows.

\begin{prop}
\label{prop: 4.3}
If $(A, J)$ is a flip pair, then
$$p_{2m-1, 0}(\sigma_A, \varphi_{J, A}) = {\Delta_J}^\textsf{T} \left( A^{m-1} \right) \Delta_{AJ},$$
$$p_{2m, 0}(\sigma_A, \varphi_{J, A}) = {\Delta_J}^\textsf{T} \left( A^m \right) \Delta_J \quad \text{and}$$
$$p_{2m, 1}(\sigma_A, \varphi_{J, A}) = {\Delta_{JA}}^\textsf{T} \left( A^{m-1} \right) \Delta_{AJ} \quad (m=1, 2, \cdots).$$
\end{prop}

\section{The Relationships between The Invariants}
\label{sec:fifth}

In Section 2, we showed that two $D_{\infty}$-TMCs $(\textsf{X}_A, \sigma_A, \varphi_{A, J})$ and $(\textsf{X}_B, \sigma_B, \varphi_{B, K})$ are $D_{\infty}$-conjugate if and only if there is a $D_{\infty}$-SSE of lag $2l$ between $(A, J)$ and $(B, K)$ for some positive integer $l$. In Section 3, we introduced a notion of $D_{\infty}$-SE, which is a conjugacy invariant of $D_{\infty}$-TMCs. The Lind zeta function discussed in Section 4 is also a conjugacy invariant. In this section, we will show that $D_{\infty}$-SE and the Lind zeta function are not complete conjugacy invariants by proving Theorem B.
The following example proves (a) of Theorem B.

\begin{exa}
\label{exa: 1}
Let $(A, I)$ and $(A, J)$ be flip pairs given by
$$A= \left[\begin{array}{rr} 1 & 1 \\ 1 & 1\end{array}\right], \qquad I=\left[\begin{array}{rr} 0 & 1 \\ 1 & 0 \end{array}\right] \qquad \text{and} \qquad J=\left[\begin{array}{rr} 0 & 1 \\ 1 & 0 \end{array}\right].$$
Direct computations yields that 
$$(A^l, A^l): (A, I) \sim (A, J) \; (\text{lag } 2l)$$ 
for all positive integers $l$.
The corresponding $D_{\infty}$-TMCs $(\textsf{X}_A, \sigma_A, \varphi_{A, I})$ and $(\textsf{X}_A, \sigma_A, \varphi_{A, J})$, however, do not share the same Lind zeta functions:
$$\zeta_{A, I}(t) = \frac{1}{\sqrt{1-2t^2}} \, \exp \left( \frac{2t+3t^2}{1-2t^2} \right),$$
$$\zeta_{A, J}(t) = \frac{1}{\sqrt{1-2t^2}} \, \exp \left( \frac{t^2}{1-2t^2} \right).$$
\end{exa}

The following example proves (b) of Theorem B.
\begin{exa}
Let
$$A=\left[\begin{array}{rrrrrrr} 1 & 1& 1 & 0 & 0 & 0 & 0 \\ 0 & 1 & 0 & 1 & 0 & 0 & 0 \\ 0 & 0 & 1 & 0 & 0 & 1 & 0 \\ 0 & 0 & 0 & 1 &0 & 0 & 1 \\ 1 & 1 & 1 & 0 & 1 & 0 & 0 \\ 1 & 1 & 1 & 0 & 0 & 1 & 0 \\ 0 & 0 & 0 & 1 & 1 & 0 & 1\end{array}\right], \quad B=\left[\begin{array}{rrrrrrr} 1 & 1 & 0 & 0 & 0 & 0 & 0 \\ 0 & 1 & 1 & 0 & 1 & 0 & 0 \\ 0 & 0 & 1 & 0 & 0 & 1 & 1 \\ 0 & 0 & 0 & 1 & 0 & 1 & 1 \\ 1 & 1 & 0 & 0 & 1 & 0 & 0 \\ 0 & 0 & 0 & 0 & 1 & 1 & 0 \\ 0 & 0 & 0 & 1 & 0 & 0 & 1 \end{array}\right],$$
$$C=\left[\begin{array}{rrrrrrr} 1 & 1& 0 & 0 & 0 & 0 & 0 \\ 0 & 1 & 0 & 1 & 1 & 1 & 0 \\ 0 & 0 & 1 & 1 & 1 & 1 & 0 \\ 0 & 0 & 0 & 1 &0 & 0 & 1 \\ 1 & 0 & 0 & 0 & 1 & 0 & 0 \\ 0 & 0 & 1 & 0 & 0 & 1 & 0 \\ 0 & 0 & 0 & 1 & 1 & 1 & 1\end{array}\right] \quad \text{and} \quad J=\left[\begin{array}{rrrrrrr} 1 & 0 & 0 & 0 & 0 & 0 & 0 \\ 0 & 0 & 0 & 0 & 1 & 0 & 0 \\ 0 & 0 & 0 & 0 & 0 & 1 & 0 \\ 0 & 0 & 0 & 0 & 0 & 0 & 1 \\ 0 & 1 & 0 & 0 & 0 & 0 & 0 \\ 0 & 0 & 1 & 0 & 0 & 0 & 0 \\ 0 & 0 & 0 & 1 & 0 & 0 & 0 \end{array}\right].$$
Then $(A, J)$, $(B, J)$ and $(C, J)$ are flip pairs.
The corresponding $D_{\infty}$-TMCs to the flip pairs share the same Lind zeta functions:
$$\sqrt{\frac{1}{t^2(1-t^2)^4(1-3t^2+t^4)}} \exp \Big(\frac{t+3t^2-t^3-2t^4}{1-3t^2+t^4}\Big).$$
Actually, they have the same numbers of fixed points:
\begin{eqnarray*}
p_m &=& 4 + \lambda^m + \mu^m,\\
p_{2m-1, 0} &=& \frac{8 \lambda^{m}-3\lambda^{m-1}}{11 \lambda -4} + \frac{8 \mu^{m}-3 \mu^{m-1}}{11 \mu -4},\\
p_{2m, 0} &=& \frac{\lambda^{m+1}}{11 \lambda -4} + \frac{\mu^{m+1}}{11 \mu -4},\\
p_{2m, 1} &=& \frac{55\lambda^m -21\lambda^{m-1}}{11 \lambda -4} + \frac{55\mu^m -21\mu^{m-1}}{11 \mu -4} \qquad (m=1, 2, \cdots).
\end{eqnarray*}
Here, $\lambda$ and $\mu$ are the zeros of $t^2 - 3t +1$:
$$\lambda = \frac{3+\sqrt{5}}{2} \qquad \text{and} \qquad \mu = \frac{3-\sqrt{5}}{2}.$$ 

If there is a $D_{\infty}$-SE $(R, S)$ between two flip pairs $(A, J)$ and $(B, K)$, then $(R, S)$ also becomes a SE between $A$ and $B$. 
Direct computations show that $A$ and $B$ are shift-equivalent. 
The matrices $A$ and $B$ have the same Jordan canonical forms up to the order of Jordan blocks while, Jordan canonical form of $C$ is different from them.
This implies that $A \nsim C$ and $B \nsim C$. (For more details, see \cite{PW} or Section 7.4 of \cite{LM}.) 
From this, we see that $(C, J)$ cannot be $D_{\infty}$-shift equivalent to $(A, J)$ or $(B, J)$.

Now, we show that if $R$ and $S$ satisfy
$$A^{2l}=RS, \quad B^{2l}=SR, \quad AR=RB \quad \text{and} \quad S=JR^{\textsf{T}}J$$
for some positive integer $l$,
then $R$ and $S$ are not integral matrices.

Let $\mathcal{J}$ denote the Jordan canonical form of $A$ or $B$: 
$$\mathcal{J} = \left[\begin{array}{rrrrrrr} \lambda & & & & & & \\ & \mu & & & & & \\ & & 0 & & & & \\ & & & 1 & 1 & 0 & 0 \\ & & & 0 & 1 &1 & 0 \\ & & & 0 & 0 & 1 &1 \\ & & & 0 & 0 & 0 & 1 \end{array}\right].$$  
The basic theorem on Jordan canonical form tells us that there are non-singular real matrices $P$ and $Q$ such that
$$A=P\mathcal{J}P^{-1} \qquad \text{and} \qquad B=Q\mathcal{J}Q^{-1}.$$
Since $AR=RB$, it follows that $A(RQ)=(RQ)\mathcal{J}$. Similarly, from $SA=BS$, it follows that $B(SP)=(SP)\mathcal{J}$.
Since $\det(A)$ and $\det(B)$ are non-zero and $A^{2l}=RS$ and $B^{2l}=SR$,
$R$ and $S$ are non-singular.
Thus, 
$$A = (RQ)\mathcal{J}(RQ)^{-1} \qquad \text{and} \qquad B = (SP)\mathcal{J}(SP)^{-1}.$$

Direct calculations yield that $RQ$ and $SP$ have the following forms:
$$RQ = P\left[\begin{array}{rrrrrrr} r_1 & 0 & 0 & 0 & 0 & 0 & 0 \\ 0 & r_2 & 0 & 0 & 0 & 0 & 0 \\ 0 & 0 & r_3 & 0 & 0 & 0 & 0 \\ 0 & 0 & 0 & r_4 & r_5 & r_6 & r_7 \\ 0 & 0 & 0 & 0 & r_4 & r_5 & r_6 \\ 0 & 0 & 0 & 0 & 0 & r_4 & r_5 \\ 0 & 0 & 0 & 0 & 0 & 0 & r_4 \end{array}\right]$$
and
$$SP = Q\left[\begin{array}{rrrrrrr} s_1 & 0 & 0 & 0 & 0 & 0 & 0 \\ 0 & s_2 & 0 & 0 & 0 & 0 & 0 \\ 0 & 0 & s_3 & 0 & 0 & 0 & 0 \\ 0 & 0 & 0 & s_4 & s_5 & s_6 & s_7 \\ 0 & 0 & 0 & 0 & s_4 & s_5 & s_6 \\ 0 & 0 & 0 & 0 & 0 & s_4 & s_5 \\ 0 & 0 & 0 & 0 & 0 & 0 & s_4 \end{array}\right]$$
for some real numbers $r_1, \cdots, r_7$ and $s_1, \cdots, s_7$.
From $S=JR^{\textsf{T}}J$, $A^{2l}=RS$ and $B^{2l}=SR$,
it follows that 
$$r_1^2 = \lambda^l, \quad r_2^2 = \mu ^l, \quad r_3^2=0,$$
$$r_4 = \pm 1, \quad r_5 = \pm l, \quad r_6 = \pm \frac{l(l-1)}{2} \quad \text{and} \quad r_7 = \pm \frac{l(l-1)(l-2)}{6}$$
and it can be shown that all the entries of $R$ and $S$ are not integers.
\end{exa}

\section{A Further Question}

The existence of SSE between two defining matrices is a necessary and sufficient condition for the corresponding shifts of finite type to be conjugate. However, there is no known algorithm for deciding whether two matrices are strong shift equivalent. Classification of shifts of finite type up to conjugacy has been one of the central problems in symbolic dynamics for many years. 

Let 
$$A= \left[\begin{array}{rrrrrrrr} 1 & 1 & 0 & 0 & 0 & 0 & 0 & 0\\ 0 & 0 & 1 & 0 & 0 & 0 & 1 & 0 \\ 0 & 0 & 0 & 1 & 0 & 1 & 0 & 0 \\ 0 & 1 & 0 & 0 & 0 & 0 & 0 & 1 \\ 1 & 0 & 0 & 0 & 1 & 0 & 0 & 0\\ 0 & 0 & 0 & 0 & 1 & 0 & 0 & 1 \\ 0 & 0 & 1 & 0 & 0 & 1 & 0 & 0 \\ 0 & 0 & 0 & 1 & 0 & 0 & 1 & 0\end{array}\right] \qquad \text{and} \qquad B=\left[\begin{array}{rr} 1 & 1 \\ 1 & 1\end{array}\right].$$

The TMC $(\textsf{X}_A, \sigma_A)$ is known as Ashley's eight-by-eight (See Problem 3.2. (2) of \cite{B} or Example 2.2.7 of \cite{Ki2}.)
and $(\textsf{X}_B, \sigma_B)$ is the full 2-shift.
The following question is due to J. Ashley in 1989: 

\begin{que}\label{que: 7}
Is $(\textsf{X}_A, \sigma_A)$ conjugate to $(\textsf{X}_B, \sigma_B)$? 
\end{que}

As in the case of shifts of finite type, it is currently unknown how to classify $D_{\infty}$-TMCs up to $D_{\infty}$-conjugacy.

If we set
$$ J= \left[\begin{array}{rrrrrrrr} 0 & 0 & 0 & 0 & 1 & 0 & 0 & 0\\ 0 & 0 & 0 & 0 & 0 & 1 & 0 & 0 \\ 0 & 0 & 0 & 0 & 0 & 0 & 1 & 0 \\ 0 & 0 & 0 & 0 & 0 & 0 & 0 & 1 \\ 1 & 0 & 0 & 0 & 0 & 0 & 0 & 0\\ 0 & 1 & 0 & 0 & 0 & 0 & 0 & 0 \\ 0 & 0 & 1 & 0 & 0 & 0 & 0 & 0 \\ 0 & 0 & 0 & 1 & 0 & 0 & 0 & 0\end{array}\right], \qquad I=\left[\begin{array}{rr} 1 & 0 \\ 0 & 1 \end{array}\right] \qquad \text{and} \qquad K=\left[\begin{array}{rr} 0 & 1 \\ 1 & 0\end{array}\right],$$
then $(A, J)$, $(B, I)$ and $(B, K)$ become flip pairs.
In Example 1, we calculated the Lind zeta functions of $(\textsf{X}_B, \sigma_B, \varphi_{B, I})$ and $(\textsf{X}_B, \sigma_B, \varphi_{B, K})$ and  showed that they are not $D_{\infty}$-conjugate. 
However, a direct calculation tells us that the Lind zeta functions of $(\textsf{X}_A, \sigma_A, \varphi_{A, J})$ and $(\textsf{X}_B, \sigma_B, \varphi_{B, K})$ coincide. Moreover, if we set
$$R = 2\left[\begin{array}{rr} 1 & 1 \\ 1 & 1 \\ 1 & 1 \\ 1 & 1 \\ 1 & 1 \\ 1 & 1 \\ 1 & 1 \\ 1 & 1 \end{array}\right]\qquad \text{and} \qquad S= 2\left[\begin{array}{rrrrrrrr} 1 & 1 & 1 & 1 & 1 & 1 & 1 & 1 \\ 1 & 1 & 1 & 1 & 1 & 1 & 1 & 1 \end{array}\right],$$
then $(R, S)$ is a $D_{\infty}$-SE of lag $6$ from $(A, J)$ to $(B, K)$.

We conclude the paper with the following question:

\begin{que}\label{que: 7}
Is $(\textsf{X}_A, \sigma_A, \varphi_{A, J})$ is $D_{\infty}$-conjugate to  $(\textsf{X}_B, \sigma_B, \varphi_{B, K})$? 
\end{que}

\end{document}